\documentclass[12pt]{article}

\usepackage{latexsym}
\usepackage{amssymb,latexsym,amsfonts,amsmath,amsthm,verbatim}
\usepackage[mathscr]{eucal}
\newtheorem{thm}{Theorem}
\newtheorem{lem}[thm]{Lemma}

\newtheorem{lemma}[thm]{Lemma}

\newtheorem{prop}[thm]{Proposition}
\newtheorem{cor}[thm]{Corollary}

\theoremstyle{remark}
\newtheorem*{rmk}{Remark}
\newcommand{\PP}{\mathbb{P}}
\newcommand{\EE}{\mathbb{E}}
\newcommand{\RR}{\mathbb{R}}

\def\eps{\varepsilon}
\def\sph{S^{n-1}}

\def\qed{\hfill $\vcenter{\hrule height .3mm
\hbox {\vrule width .3mm height 2.1mm \kern 2mm \vrule width .3mm
height 2.1mm} \hrule height .3mm}$ \bigskip}
\def\P{\mathbb{P}}

\def\R{\mathbb{R}}
\def\PP{\mathbb{P}}
\def\EE{\mathbb{E}}
\def\RR{\mathbb{R}}

\def\Tr{\operatorname{Tr}}

\begin{document}

\title{Polynomial bounds \\ for large Bernoulli sections of $\ell_1^N$}
\author{S. Artstein-Avidan\textsuperscript{1,2},
        O. Friedland\textsuperscript{3},
        V. Milman\textsuperscript{1,3},
        S. Sodin\textsuperscript{3}}
\footnotetext[1]{partially supported by BSF grant 2002-006}
\footnotetext[2]{supported by the National Science Foundation under agreement No. DMS-0111298}
\footnotetext[3]{supported in part by the Israel Science Academy}
\date{}
\maketitle

\section{Introduction}

This paper consists of two distinct parts. The first one presents
the ``local'' version of the result of Bai and Yin from \cite{BY}.
This result gives an estimate from below for the probability that
the smallest singular value of a random sign matrix is outside some
interval. In particular, it gives a lower bound for the probability
that an ``almost square'' matrix, that is, a $(1-\delta)n\times n$
matrix, has smallest singular value above $\approx \delta$. This is
a ``finite dimensional'' version of the results of Bai and Yin
\cite{BY}, and in this ``local'' version it is much more useful for
applications in Asymptotic Geometric Analysis problems, where
quantitative estimates of deviations are needed. This is presented
in Section~\ref{sash}. A more extensive presentation of this result
will be given by the IVth named author in \cite{safut}.

The second part of this paper consists of precisely such an application,
where the method of \cite{AFM} and some other recent developments are joined
with the above, to improve results from \cite{LPRTV} and from \cite{AFM}
regarding the distance from euclidean space of almost full
dimensional sections of the space $\ell_1^N$ realized as images of
sign matrices. For $N= (1+\delta)n$ we receive estimates on the
isomorphism constant which are much better than were previously
known, and in particular are {\em polynomial} in $\delta$.

{\em Acknowledgments} The authors thank Prof.~S.~Szarek for bringing the
paper \cite{BY} to their attention and suggesting that the methods used there
yield an estimate on the probability of deviation. The IVth named author
thanks Prof.~E.~Gluskin for useful discussions of the combinatorial estimates
in Section~\ref{gt}.
The authors thank the referee for his careful reading and useful comments.

\section{The rate of convergence in the result of Bai and Yin}\label{sash}

In this section we present a lower bound on the least singular value
of a Bernoulli random matrix, in the spirit of Bai and Yin \cite{BY}.

\subsection{Introduction and main statement}

Let $X$ be a $p \times n$ matrix of random signs: $X_{ik}$ are
independent for $1 \leq i \leq p$ and $1 \leq k \leq n$,
\begin{equation}\label{def_x}
\PP \left\{ X_{ik} = 1\right\}
    = \PP \left\{ X_{ik} = -1 \right\}
    = 1/2 \, \text{.}
\end{equation}

We study the spectrum $\Lambda_{S}$ of the covariance matrix
\begin{equation}\label{def_s}
S = n^{-1} X X^\text{T} \, \text{.}
\end{equation}

Let $\mu_{S} = p^{-1} \sum_{\lambda \in \Lambda(S)} \delta_\lambda$
be the empirical eigenvalue distribution of $S$. Marchenko and
Pastur \cite{mp} proved that
\[ d \mu_S \overset{\text{a.s.}}{\longrightarrow} f_\text{MP} dx
    \quad \text{as} \quad n \to \infty \, \text{,}\]
where the limit density equals
\begin{equation}\label{mp}
f_{\text{MP}}(x) =
    \begin{cases}
    \frac{1}{2 \pi y x}
        \sqrt{ (x - a) (b - x) },
        &a \leq x \leq b \\
    0, &\text{otherwise,}
    \end{cases}
\end{equation}
%  I added a different word this time, is it OK?
%% Tayib
with the notation
\begin{equation}\label{def_yab}
y = p/n < 1 \:\:\text{(fixed)}, \quad a = (1 - \sqrt{y})^2, \quad b = (1 + \sqrt{y})^2 \,
\text{.}
\end{equation}

It is natural to ask whether the eigenvalues of $S$ can lie far from
the support $[a, \, b]$ of this distribution. Bai and Yin \cite{BY}
answered negatively, proving (for a more general random matrix
model) that with probability 1
\[ \lambda_{\min}(S) \to a, \quad \lambda_{\max}(S) \to b
    \quad \text{as $n \to \infty$.}  \tag{BY} \]

In the spirit of local theory we strive for a quantitative form of
this result.

\begin{thm}\label{thm}
There exists a universal constant $C > 0$ such that the following
holds. Let $X$ be a $p \times n$ matrix of random signs as defined
by (\ref{def_x}); define $S$ as in (\ref{def_s}) and $y$, $a$ and $b$
as in (\ref{def_yab}); assume that
\begin{equation}\label{restr_thm}
\frac{C\log^2 n}{\sqrt{y}\sqrt[3]{n}} \leq \epsilon \leq 1 \, \text{.}
\end{equation}
Then the probability that $S$ has eigenvalues outside $[a - \epsilon, \,
b + \epsilon]$ is less than
\[ \exp \left( - C^{-1} y^{1/6} n^{1/6} \epsilon^{1/2} \right)
    = \exp \left( - C^{-1} p^{1/6} \epsilon^{1/2} \right)
    \, \text{.} \]
\end{thm}

For $y$ close to 1 the theorem yields the following lower bound on
the least eigenvalue of $S$:
\begin{thm}\label{res} There exists a universal constant $C > 0$
such that if, in the notation of Theorem \ref{thm},
$y = 1 - \delta$ with $1/2 > \delta > C n^{-1/6} \log n$, then
\[ \PP \left\{ \lambda_{\min}(S) \leq \delta^2/8 \right\}
    \leq \exp \left(- C^{-1} n^{1/6} \delta \right) \, \text{.}\]
\end{thm}

Recently, Litvak, Pajor, Rudelson and Tomczak-Jaegermann
\cite{lprt} proved (in a more general setting) that
if $y = 1 - \delta$ with $1>\delta \geq c_1/\ln c_2n$ in the notation
of Theorem \ref{thm}, then
\[ \PP \left\{ \lambda_\text{min}(S) \leq A a^{1/\delta} \right\}
    \leq \exp(-Cn)
    \, \text{,} \label{LPRT} \tag{LPRT} \]
where $A>0$, $1>a>0$, $C,c_1,c_2>0$ are universal constants.

Note that the bound on the probability decays exponentially; this is
rather important in geometric applications. We do not know whether
the left-hand side in Theorem~\ref{res} is in fact as small as
$\exp \left(- n\delta^C/C \right)$ for some $C > 0$.

Let us show that Theorem \ref{thm} implies Theorem \ref{res}.

\begin{proof}[Proof of Theorem \ref{res}]
The Taylor expansion yields $\sqrt{y} \approx 1 - \delta/2$ and
hence
\[ (1-\sqrt{y})^2 - \epsilon \approx \delta^2/4 - \epsilon \, \text{.} \]
Now take $\epsilon \approx \delta^2/8$ and use Theorem \ref{thm}. We
obtain:
\[ \PP \left\{ \lambda_\text{min}(S) < \delta^2/8 \right\}
    \leq \exp \left( - \frac{\sqrt[6]{y}}{\sqrt{8}C} \, n^{1/6} \delta \right)
    \leq \exp \left( - \frac{n^{1/6} \delta}{2^{5/3} C} \right)
    \, \text{.} \]
\end{proof}

The main idea behind the proof of Theorem \ref{thm} makes use of the
following construction, due to Bai and Yin \cite{BY}. We define a
sequence of matrices $T(l)$, $l = 0, 1, 2, \cdots$, that are certain
polynomials of the matrix $T = S - \mathbb{I}$: $T(l) = p_l(T)$. If
$\mu_1, \cdots, \mu_p$ are the eigenvalues of $T$, then $p_l(\mu_1),
\cdots, p_l(\mu_p)$ are the eigenvalues of $T(l)$.

The polynomials $p_l$ can be expressed via the Chebyshev polynomials
of the second kind. If $\mu \notin [a-1, \, b-1]$, the sequence
$p_l(\mu)$ tends to infinity exponentially fast. We define $p_l$ and
prove these observations in Section~\ref{cheb}.

On the other hand, the expression $\EE \Tr T(l)$ allows a
graph-theoretical interpretation showing that it can not grow too
fast. We prove such a bound in Section~\ref{gt}, using a
modification of the combinatorial argument due to Bai and Yin.

In Section~\ref{endofproof} we combine these facts and obtain a bound
on $a - \lambda_{\min}(S)$, $\lambda_{\max}(S) - b$ that concludes
the proof of Theorem~\ref{thm}.

\subsection{Construction and basic properties of $T(l)$}\label{cheb}

Denote $y_1 = \frac{p-2}{n}$, $y_2 = \frac{(p-1)(n-1)}{n^2}$;
$y \geq y_2 \geq y_1 = y - \frac{2}{n}$.

Define a sequence of matrices $T(l) = (T_{ij}(l))_{ij}$,
\begin{equation}\begin{cases}\label{tdef}
&T(0) = \mathbb{I}, \quad T(1) = T = n^{-1}X X^\text{T} - \mathbb{I}, \\
&T(l+1) =  \left(T - y_1\mathbb{I} \right) \cdot T(l)
                - y_2\cdot T(l-1)\,\text{.}
\end{cases}\end{equation}

We have: $T(l) = p_l(T)$, where
\[\begin{cases}
&p_0(\mu) = 1, \quad p_1(\mu) = \mu, \\
&p_{l+1}(\mu) = \left( \mu - y_1 \right) \cdot p_l(\mu)
                - y_2 \cdot p_{l-1}(\mu) \, \text{.}
\end{cases}\]

Recall the definition
\[ U_l(\cos \theta) =
    \frac{\sin \left((l+1) \, \theta \right)}{\sin \theta}
    \tag{Cheb1}\label{Cheb1}\]
of the Chebyshev polynomials of the second kind. Here, both the right-hand
side and the left-hand side are polynomials; hence the equality makes sense
for any $\theta \in \mathbb{C}$.

Equivalently, $U_l$ can be defined by the recurrence relation
\[\begin{cases}
&U_0(x) = 1, \quad U_1(x) = 2x, \\
&U_{l+1}(x) = 2x U_l(x) - U_{l-1}(x) \, \text{.}
\end{cases}\tag{Cheb2}\label{Cheb2}\]

The latter definition readily yields the formula
\begin{equation}
p_l(\mu) = {y_2}^{l/2} U_l\left(\frac{\mu - y_1}{2 \sqrt{y_2}}\right) +
          y_1 {y_2}^{(l-1)/2} U_{l-1}\left(\frac{\mu - y_1}{2 \sqrt{y_2}}\right)
          \, \text{.}
\label{viacheb}\end{equation}

\begin{rmk}
If we replace $y_1$ and $y_2$ with $y$ in (\ref{viacheb}), the
sequence becomes orthogonal with respect to the Marchenko--Pastur
measure (\ref{mp}). Kusalik, Mingo and Speicher \cite{KMS} used a
different form of this sequence to study the spectral properties
of random matrices with complex Gaussian entries, and called it
the sequence of {\em shifted Chebyshev polynomials of the second
kind}.
\end{rmk}

Now we use (\ref{Cheb1}) to estimate the polynomials $p_l$.

\begin{lemma} There exists a universal constant $C>0$ such that
the following properties hold for any even $l \geq 2$, $0 \leq
\epsilon \leq 1$:
\begin{enumerate}
\item For any $\mu \in \RR$,
\begin{equation}\label{ch1}
p_l(\mu) \geq - 2 l y^{l/2} \, \text{.} \end{equation}
\item If
\[ | \mu - y_1| \geq 2\sqrt{y_2} + \epsilon \, \text{,} \]
then
\begin{equation}\label{ch2}
p_l(\mu) \geq {y_2}^{l/2}
    \exp \left( C^{-1} l \epsilon^{1/2} y_2^{-1/4} \right)
\, \text{.} \end{equation}
\end{enumerate}
\end{lemma}

\begin{proof}\hfill
\begin{enumerate}
\item If $x$ lies outside the interval
$[\cos \frac{l\pi}{l+1}, \cos \frac{\pi}{l+1}]$, then $U_l(x) > 0$.
Otherwise, $x = \cos \theta$ for some $\frac{\pi}{l+1} \leq \theta
\leq \frac{l\pi}{l+1}$; therefore
\[ U_l(x) \geq - \sin^{-1} \frac{\pi}{l+1} \geq - \frac{l+1}{2} \, \text{.} \]
Hence
\[ p_l(\mu)
    \geq - \left({y_2}^{l/2} + y_1 {y_2}^{(l-1)/2}\right) \, \frac{l+1}{2}
    \geq - 2 y^{l/2} l \, \text{.} \]
\item If $|x| \geq 1 + \epsilon$, $x = \cos i\theta$ for some
$\theta \geq C^{-1}\epsilon^{1/2}$; hence
\[ U_l(x) = \frac{\sin ((l+1)\,i\theta)} {\sin i\theta}
    \geq e^{l \theta/2} \geq e^{C_1^{-1} l \epsilon^{1/2}} \, \text{.} \]
Therefore
\[ p_l(\mu)
    \geq {y_2}^{l/2} \exp \left( \,
\frac{l}{C} \frac{\sqrt{\epsilon}}{\sqrt[4]{y}}
    \, \right) \, \text{.}\]
\end{enumerate}

\end{proof}

Next we apply (\ref{ch1}) and (\ref{ch2}) to the eigenvalues of $T$.

\begin{lemma}\label{TRGEQ} There exists a universal constant $C>0$ such that
if $n \geq l \geq 2$ is even,
\begin{equation}\label{restr_lemma}
C \max \left( \frac{1}{\sqrt{y}n} \, ,
    \frac{\sqrt{y} \log^2 n}{l^2} \right) \leq \epsilon \leq 1 \qquad \text{and}
\end{equation}
\[ \max \left\{ | \mu - y|  \: \big{|} \:
    \text{$\mu$ is an eigenvalue of $T$} \right\}
    \geq 2\sqrt{y} + \epsilon \, \text{,} \]
then
\begin{equation}
\Tr T(l)
    \geq y^{l/2} \exp \left( C^{-1} l \epsilon^{1/2} y^{-1/4} \right)
    \, \text{.}
\label{trgeq}\end{equation}
\end{lemma}

\begin{proof}
Let $\mu_1, \cdots, \mu_p$ be the eigenvalues of $T$; suppose
$|\mu_1 - y| \geq 2 \sqrt{y} + \epsilon$.

Then by (\ref{restr_lemma})
\[\begin{split}
| \mu_1 - y_1 |
    &\geq 2 \sqrt{y_2} + \epsilon - 2 \left[ \sqrt{y} - \sqrt{y_2} \right] - [ y - y_1 ] \\
    &\geq 2 \sqrt{y_2} + \epsilon - \frac{4}{n\sqrt{y}} - \frac{2}{n} \\
%%% didn't check it yet
%% This is just (a-b)(a+b) = (a^2-b^2)
    &\geq 2 \sqrt{y_2} + \epsilon - \frac{6}{n\sqrt{y}}
    \geq 2 \sqrt{y_2} + C_1 \epsilon \, \text{.}
\end{split}\]

Write the bound (\ref{ch2}) with $\mu = \mu_1$ and the bound
(\ref{ch1}) with $\mu = \mu_2, \cdots, \mu_p$; add the
inequalities and use (\ref{restr_lemma}) once again:
%%% i still should check, with condition on epsilon and everything.
\[\begin{split}
\Tr T(l)
    &\geq {y_2}^{l/2} \exp \left( C^{-1} l
        \frac{\sqrt{C_1\epsilon}}{\sqrt[4]{y_2}} \right)
        - 2 l p \, y^{l/2} \\
    &\geq C_2^{-1} y^{l/2} \exp \left( C_2^{-1} l
        \frac{\sqrt{\epsilon}}{\sqrt[4]{y}} \right)
        - 2 n^2 y^{l/2} \\
    &\geq y^{l/2} \exp \left( C_3^{-1} l \epsilon^{1/2} y^{-1/4} \right)
        \, \text{.}
\end{split}\]
\end{proof}

\subsection{Combinatorial description of $T(l)$}\label{gt}

Now we give a combinatorial description of $\EE \Tr T(l)$.

\begin{lemma}
The following equality holds:
\begin{equation}
T_{ij}(l) =
    \frac{1}{n^l} \, \sum\nolimits^\ast X_{iv_1} X_{u_1v_1} X_{u_1v_2} X_{u_2v_2}
        \cdots X_{u_{l-1}v_l} X_{jv_l} \, \text{,}
\label{comb}\end{equation} where the sum $\sum^\ast$ is over all
$u_1, \ldots, u_{l-1}$ and $v_1, \ldots, v_{l}$ satisfying
$1 \leq u_r \leq p$ for $1 \leq r \leq l-1$ and  $1 \leq v_s \leq n$ for $1
\leq s \leq l$, and such that, in addition,
\[\begin{cases}
i \neq u_1 \neq u_2 \neq u_3 \neq \cdots \neq u_{l-1} \neq j&\\
v_1 \neq v_2 \neq v_3 \neq \cdots \neq v_l \, &\text{.}
\end{cases}\]
\end{lemma}
\noindent (Notice that there is no requirement $u_1 \neq u_3$, for example.)

\begin{proof}
Denote by $T_{ij}'(l)$ the right-hand side of (\ref{comb}); denote
$T'(l) = (T'_{ij}(l))$. Then $T'(0) = \mathbb{I} = T(0)$, $T'(1) = T
= T(1)$.

Further, $(T \cdot T'(l-1))_{ij}$ is a sum of the same form as
(\ref{comb}), but without the condition $v_1 \neq v_2$. The three
cases ({\em i})~$v_1 \neq v_2$, ({\em ii})~$v_1 = v_2$ and $i \neq
u_2$, ({\em iii})~$v_1 = v_2$ and $i = u_2$ yield the terms
\[ T'(l)\, , \quad y_1\, T'(l-1)\, , \quad y_2 \, T'(l-2)\, , \]
%%% Still keeping the  note that I only checked
%%% the middle one and not the third;
respectively. Therefore $T'(l)$ satisfy the same recurrence relation
(\ref{tdef}) as $T(l)$; this concludes the proof.
\end{proof}

The random variables $X_{uv}$ are independent; therefore the expectation
of a term in (\ref{comb}) vanishes unless every $X_{uv}$ appears an even
number of times in the product. In the latter case, the expectation
equals $1$ (note that $0$ is even).

\begin{cor}\label{cor1}
The expectation $n^l \, \EE\Tr T(l)$ equals the number of configurations
\[ 1 \leq i,u_1,u_2,u_3,\cdots,u_{l-1} \leq p, \quad
   1 \leq v_1,v_2,\cdots,v_l \leq n \, \text{,} \]
such that
\[\begin{cases}
i \neq u_1 \neq u_2 \neq u_3 \neq \cdots \neq u_{l-1} \neq i&\\
v_1 \neq v_2 \neq v_3 \neq \cdots \neq v_l \, &\text{,}
\end{cases}\] and every pair $uv$
appears an even number of times in the sequence
\[ iv_1, u_1v_1, u_1v_2, u_2v_2, \cdots, u_{l-1}v_l, iv_l \, \text{.}\]
\end{cor}

The following graph-theoretical interpretation will be of use. Every
configuration of $i$, $u_r$ and $v_s$ which is permitted
in Corollary \ref{cor1} corresponds to a closed path
$W$ in the bipartite graph $K_{p, \,n}$ such that
\begin{enumerate}
\item[(W1)]\label{C1} the path $W$ passes through every edge an even number
    of times;
\item[(W2)]\label{C2} $W$ never passes through an edge 2 times consequently
    (i.e.\ the pattern $w \to w' \to w$ is not allowed).
\end{enumerate}

\noindent (Moreover, every path begins on the left side of the graph, but we ignore
this in our estimates.)

Let $W$ be a closed path on an arbitrary graph $G$ so that (W1) and
(W2) hold. Consider $W$ as a set of triples $(w_1, w_2, r)$, where
$1 \leq r \leq 2l$, meaning that the $r$th edge on $W$ goes from $w_1$ to
$w_2$.

Divide the edges into 3 classes:
\[\begin{split}
T_1 &= \big\{ (w_1, w_2, r) \in W \, \big{|} \, \\
    &\qquad\qquad\forall r' < r, (w_1', w_2', r') \in W \Rightarrow
                    w_1' \neq w_2 \wedge w_2' \neq w_2 \big\}
                    \, \text{,}\\
T_2 &= \big\{ (w_1, w_2, r) \in W \, \big{|} \, \\
    &\qquad\qquad\exists r' < r:
        (w_1, w_2, r') \in T_1 \vee (w_2, w_1, r') \in T_1 \, \text{,} \\
    &\qquad\qquad\forall r' < r'' < r:
        (w_1, w_2, r'') \notin W \wedge(w_2, w_1, r'') \notin W \big\} \, \text{,}\\
T_3 &= W \, \backslash \, (T_1 \cup T_2) \, \text{.}
\end{split}\]

({\em Semiformal verbal description:} The edges of $T_1$ are the
first edges to visit their endpoints; that is, $T_1$ is the DFS
tree of $W$. Every edge in $T_1$ appears at least once again on
$W$; we denote by $T_2$ the set of second appearances of the $T_1$
edges. All the other edges form the set $T_3$.)

Let us call a sequence of vertices $f = w_1 w_2 \cdots w_k$ ($k > 1$)
a {\em protofragment} of $W$ if the following 3 conditions hold: {\em (i)}
for some $r$
\[ (w_1, w_2, r), (w_2, w_3, r+1), \cdots,
    (w_{k-1}, w_k, r+k-2) \in T_1 \, \text{,} \]
{\em (ii)} for some $r'$
\[\begin{cases}
\text{either}
    &(w_1, w_2, r'), \, (w_2, w_3, r'+1), \, \cdots, \, (w_{k-1}, w_k, r'+k-2) \in T_2 \\
\text{or}
    &(w_k, w_{k-1}, r'), \, \cdots, \, (w_3, w_2, r'+k-3), \,
        (w_2, w_1, r'+k-2) \in T_2 \, \text{,}
\end{cases}\]
and {\em (iii)} $f$ is maximal with respect to the 2 conditions {\em (i)--(ii)}.

If $f = w_1 w_2 \cdots w_k$ is a protofragment, $w_1 \neq i$, we
call its suffix $\bar{f} = w_2 \cdots w_k$ a {\em fragment} of
length $k-1$. If $w_1= i$, we call $f$ a {\em fragment} of length $k$.
%why did I write (k+1) before??
The vertices of $W$ are thereby divided into $F$ fragments.

The following combinatorial bound will be crucial ($\sharp$ denotes cardinality):

\begin{lemma}\label{mcl}
$F \leq 2 \sharp T_3 + 1$.
\end{lemma}

\begin{proof}

Let $f$ be a protofragment that starts with $w \neq i$; consider $2$ cases.
If $f$ is passed in the same direction in $T_1$ and $T_2$, the edge adjacent
to $w$ in one of the 2 passages is in $T_3$.

Otherwise, the last edge before the second appearance of $f$ is in $T_3$.

Let$e$ be the $T_3$ edge in either case. The map $f \mapsto e$ is at most 2--1;
hence $F-1 \leq 2 \sharp T_3$.

\end{proof}

\begin{lemma}\label{nfrag}
The number of different fragments of length $k$ in $K_{p,\,n}$ is
bounded by $2y^{-1/2} (pn)^{k/2}$.
\end{lemma}

\begin{proof}
First decide to which side of the graph does the first vertex
belong. Then choose all the vertices.
\end{proof}

Now we can bound the number of paths satisfying (W1)--(W2) on $K_{p,
\,n}$. Let $V$ be the number of (distinct) vertices on $W$.

First, choose the lengths of the fragments. This can be done in
$\binom{V}{F-1} \leq V^F / F!$ ways. Next, choose the fragments
themselves; by Lemma \ref{nfrag} this can be done in at most
$(y/4)^{-F/2} (pn)^{V/2}$ ways.

We can change the directions of the fragments in $T_2$, in $2^F$
ways. Now that the fragments are ready, glue them onto the path;
this can be done in $(2l - 2V + 1)^{2F}$ ways (just pick a place for every fragment).

Now there are $2l-2V$ vertices left. Every one of them coincides with one of the $V$
vertices that we already have. Once we choose one of the $V^{2l-2V}$ arrangements,
our path is ready.

%%% (Is the explanation any better now?)
%%% Its hard for me now to judge, I will look at it again in a day or two

Multiplying all these numbers, we see that the number $\mathcal{P}$
of paths is bounded by

\[\begin{split}
\mathcal{P} &\leq \sum_{V=1}^l \sum_{F=1}^{l}
    \frac{V^F}{F!} (y/4)^{-F/2} (pn)^{V/2}
    2^F (2l - 2V + 1)^{2F} V^{2l-2V} \\
&\leq \sum_{V=1}^l \sum_{F=1}^{l} \:
    (pn)^{V/2} \left(CVy^{-1/2}\right)^F V^{2l-2V}
    \times \left( \frac{2l-2V+1}{F} \right)^F
    \, \text{.}
\end{split}\]

Now, $(x/F)^F \leq e^{x}$;
% Wouldn't it make more sense just $(x/F)^F \leq e^{x}$?
%% This does not matter at all ($e^{100 x}$ would be as good);
%% but do you know any way to prove the weaker form without
%% using the stronger one?
$F \leq 2 \sharp T_3 + 1 = 4l - 4V +5$.
Therefore

\[\begin{split}
\mathcal{P} &\leq \sum_{V=1}^l \sum_{F=1}^{l} \:
    (pn)^{V/2} (C' V y^{-1/2})^{4l - 4V + 5} V^{2l-2V} \\
&\leq \sum_{V=1}^l l \:
    (pn)^{V/2} (C' V^{3/2} y^{-1/2})^{4l - 4V + 5} \\
&\leq l^{9} y^{-5/2} (pn)^{l/2} \sum_{V=1}^l \:
    (pn)^{(V-l)/2} (C' V^{3/2} y^{-1/2})^{4l - 4V} \, \text{.}
\end{split}\]

Now, if $(C' l^{3/2} y^{-1/2})^8 < pn$, every term in the sum is no greater
than $1$. Therefore if
\[ l < C''^{-1} y^{1/3} (pn)^{1/12} = C''^{-1} y^{5/12} n^{1/6} \, \text{,} \]
then
\[ \mathcal{P} \leq l^{10} y^{-5/2} (pn)^{l/2} \, \text{;} \]
finally (in one line):
\begin{equation}
\EE \Tr T(l) \leq l^{10} y^{(l-5)/2} \quad
     \text{if $l \leq {C''}^{-1} y^{5/12} n^{1/6}$.}
     \label{trleq}
\end{equation}

\subsection{Conclusion of the proof}\label{endofproof}

\begin{proof}[Proof of Theorem \ref{thm}]
Let $l = 2  \lfloor (2C'')^{-1} y^{5/12} n^{1/6} \rfloor$ in
(\ref{trleq}).

Then by (\ref{restr_thm})
\[ \epsilon \geq \frac{C \log^2 n}{\sqrt{y} \sqrt[3]{n}}
    \geq \frac{C}{\sqrt{y}n} \]
and
\[ \epsilon \geq \frac{C \log^2 n}{\sqrt{y} \sqrt[3]{n}}
    \geq \frac{C \sqrt{y} \log^2 n}{l^2} \frac{l^2}{y\sqrt[3]{n}}
    \geq \frac{C_1 \sqrt{y} \log^2 n}{l^2} \, \text{;} \]
therefore (\ref{restr_lemma}) holds.

By Lemma~\ref{TRGEQ}, Chebyshev's inequality, the estimate (\ref{trleq})
and the condition (\ref{restr_thm}) that we imposed on $\epsilon$,
\[\begin{split}
&\PP \left\{ \text{$S$ has eigenvalues outside
    $[a - \epsilon, \, b + \epsilon]$}\right\} \\
&\qquad= \PP \left\{ \text{$T$ has eigenvalues outside
    $[a - 1 - \epsilon, \, b - 1 + \epsilon]$}\right\} \\
&\qquad\leq \PP \left\{ \Tr T(l) \geq
    y^{l/2} \exp \left( C^{-1} l \epsilon^{1/2} y^{-1/4} \right) \right\} \\
&\qquad\leq \frac{\EE \Tr T(l)}
                 {y^{l/2} \exp \left( C^{-1} l \epsilon^{1/2} y^{-1/4} \right)}
\leq C y^{5/3} n^{5/3}
    \exp \left( - C^{-1} y^{1/6} n^{1/6} \epsilon^{1/2} \right) \\
&\qquad\leq \exp \left( - {C_1}^{-1} n^{1/6} y^{1/6} \epsilon^{1/2} \right)
= \exp \left( - {C_1 }^{-1} p^{1/6} \epsilon^{1/2}
    \right) \, \text{.}
\end{split}\]

We are done.

\end{proof}

\section{Application to large sections of $\ell_1^N$}\label{afm}
It is well known that $\ell_1^{(1+\delta)n}$ has isomorphic
euclidean sections of dimension $n$ (see \cite{Kashin}), with
constant of isomorphism independent of the dimension $n$ and depending
only on $\delta$. When the section is taken to be the image of a
matrix with i.i.d. gaussian entries (which is the same as taking a
random subspace in the Grassmanian $G_{N,n}$ with respect to the
normalized Haar measure), the dependence is polynomial in $\delta$,
with high probability on the choice of the entries. This was
discovered first in the results of \cite{GG}.

The image of a matrix of signs is simply the span of some set of
vertices of the unit cube, and thus has more structure, and is
sometimes more useful in implementations. Schechtman showed in
\cite{SchF} that the image of a matrix whose rows are $N=Cn$
sign-vectors in $\R^n$, where $C$ is a universal constant, also
realizes, with high probability on the choice of signs, an isomorphic
to euclidean section of $\ell_1^N$. The question then remained
whether the constant $C$ can be reduced to be close to $1$. This was
resolved by Johnson and Schechtman, and follows from their paper
\cite{JS}. However, they showed the {\em existence} of such a sign
matrix, and not that it is satisfied for a matrix whose rows are $N
= (1+\delta)n$ {\em random} sign-vectors. In a paper by Litvak,
Pajor, Rudelson, Tomczak-Jaegermann and Vershynin \cite{LPRTV} this
was demonstrated. However, the dependence of the constant of
isomorphism on $\delta$ in their result is exponentially bad, and
they get $c(\delta) = c^{1/\delta}$. In this paper we get a better
dependence, polynomial in $\delta$, however the probability that we
get is slightly smaller than the probability in \cite{LPRTV}, with
$n^{1/6}$ in the exponent instead of $n$.

%Here if Sasha's result improves.. change...

We remark that results of this type can be viewed also in a
different way, as a realization of Khinchine inequality with few
vectors. The classical Khinchine inequality states that (for best
constants as below see \cite{Sza})
\[ {1\over {\sqrt 2}}
(\sum_{i=1}^n x_i^2)^{1/2}\le Ave_{\eps_1,\ldots,\eps_{n}~={\pm
1}}|\sum_{i=1}^n\eps_i x_i| \le (\sum_{i=1}^n x_i^2)^{1/2}.
\]
Instead of averaging over {\em all} sign-vectors we may average over
only $n(1+\delta)$ of them (chosen randomly, and good for {\em all}
$x$), and get the same inequality only with a worse constant instead
of $\sqrt{2}$. The constant is universal for fixed $\delta$, and the
way it behaves when $\delta \to 0$ is the subject of this paper,
reformulated.

\vskip 6pt

In this section we show that for a random $N \times n$ sign matrix,
where $N = n(1+\delta)$, we have with high probability that the
section of $\ell_1^N$ given by its image is isomorphic to the
euclidean ball with polynomial dependence of the constants of
isomorphism on $\delta$. The developments which allowed this
advancement include the methods of Schechtman to get $L_1$ splitting
as in \cite{split}, the quantitative version of the result of Bai
and Yin \cite{BY} given in Theorem \ref{res} of the previous
section, and the use of Chernoff bounds for geometric purposes much
like is done in \cite{AFM}. We prove

\begin{thm}\label{nice}
There exist universal constants $\delta_0$, $c',c''$, and $c_0$ such
that the following holds. Let $c''
n^{-1/6}\log n <\delta<\delta_0$, and denote $N = (1+ \delta) n$. Then with
probability greater than $1- e^{-c'\delta n^{1/6}}$, for
$(1+\delta)n$ random sign-vectors $\eps_j \in \{-1, 1\}^n$, $j = 1,
\ldots , n+\delta n$, one has for every  $x \in \R^{n}$
\begin{equation}\label{good}
c(\delta) |x| \le {1\over N}\sum_{j=1}^N | \langle x, \eps_j \rangle
|
%\le C |x|
,\end{equation} where $c(\delta) = c_0\delta^{5/2}/\log(1/\delta)$.
\end{thm}

\noindent In fact it is easy to see that once we know Theorem
\ref{nice} the above remains true for any $\delta>0$, and the
restriction $\delta<\delta_0$ is artificial. Also, an upper bound in
(\ref{good}) is known and standard, similar to Lemma \ref{uptwo}.

Notation: We pick the $N = n+\delta n$ random sign vectors $\eps_j$,
normalize them to be unit vectors by dividing by $\sqrt{n}$ and
denote the normalized vectors by $v_1, \ldots, v_{n+\delta n/2}$,
$w_1, \ldots, w_{\delta n/2}$, that is, $v_j = \eps_j/\sqrt{n}$ for
$j = 1, \ldots, n+\delta n/2$ and $w_j = \eps_{(n+\delta
n/2+j)}/\sqrt{n}$ for $j= 1, \ldots, \delta n/2$. Throughout the
proof $c$, $c_1$, $c_2'$, $C_3$ etc. will denote universal constants
which can be easily estimated.

\vskip 6pt

Our proof mimics the proof of the theorem when the first $n$ vectors
form an orthonormal basis, and then the upper square in the matrix
is an isometry. To substitute this fact, we will first of all need
an estimate for the smallest eigenvalue of an $n\times
(1+\delta/2)n$ matrix of random signs, which is given in Proposition
\ref{sasha} below, which is simply a reformulation of Theorem
\ref{res}. It can be looked upon as a near-orthogonality result for
the $n$ random column vectors which are sign-vectors that live in
$(n+\delta n/2)-$dimensional space.

\begin{prop}\label{sasha}
There exist universal constants $\delta_0$, $c'', c_1$ and $c_1'$
such that for any $c'' n^{-1/6} \log n <\delta<\delta_0$,  if $v_j$
are $n+\delta n/2$ random vectors chosen uniformly and independently
in $\{-1/\sqrt{n}, 1/\sqrt{n}\}^n$ then with probability greater
than $1- e^{-c_1'\delta n^{1/6}}$ we have for every $x \in \R^n$
that
\[
c_1\delta |x| \le \left( \sum_{j=1}^{n+\delta n} | \langle x, v_j
\rangle |^2 \right)^{1/2}.
\]
\end{prop}

%We will also need the following upper bound, which is well known and
%follows from Bernstein's inequality and successive approximation.

%\begin{lem}\label{upone} There exist universal constants $C_3$ and
%$c_3'$ such that for any $\delta
%> 0$ if $v_j$ are $n+\delta n$ random vectors of $\pm 1/\sqrt{n}$
%then with probability greater than $1- e^{-c_3'n}$ we have for every
%$x \in \R^n$ that
%\begin{equation}\label{fl}
%\left(\sum_{j=1}^{n+\delta n} | \langle x, v_j \rangle
%|^2\right)^{1/2}\le C_3 |x|.
%\end{equation}
%\end{lem}

%\vskip 6pt

The idea of the proof of Theorem \ref{nice} is to use the ``near
orthogonality'' of the first $n+\delta n/2$ row vectors to ensure a
lower bound in most directions. For the directions which remain, we
obtain a lower bound by using the last $\delta n/2$ rows. To this
end we will use a net argument, and hence we also need an upper
bound for the contribution of the last $\delta n/2$ rows. This is
given by the following

\begin{lem}\label{uptwo}
There exist universal constants $c_3'$ and $C_3$ such that for any
$\delta
> 0$ if $w_j$ are $\delta n/2$ random vectors of $\pm 1/\sqrt{n}$
then with probability greater than $1- e^{-c_3'n}$ we have for every
$x \in \R^n$ that
\begin{equation}\label{fl2} {1\over \sqrt{n}} \sum_{j=1}^{\delta n/2}
| \langle x, w_j \rangle | \le C_3 \sqrt{\delta}|x|.
\end{equation}
\end{lem}
\noindent (Notice that although for a single point, in expectation,
we have (\ref{fl2}) with $\delta$ instead of $\sqrt{\delta}$, for
the probability to suffice for the whole net we need to allow
deviation of order $\sqrt{\delta}$ from the expectation.)

\noindent{\em Proof~} Bernstein inequality implies that for any
$t>1$
\[ \P [{2\over \delta n} \sum_{j=1}^{\delta n/2}|\langle x, w_j \rangle | \ge
t{|x|\over \sqrt{n}}]\le e^{-ct^2 \delta n } \]
for a universal $c$.
We pick a $1/2$-net on the sphere $\sph$ with cardinality $5^n$ and
pick $t = \sqrt{2\ln 5\over c\delta}$. Then with probability greater
than $1- 5^{-n}$ we have that for every element $x$ in the net
\[ {1\over \sqrt{n}} \sum_{j=1}^{\delta n/2}
| \langle x, w_j \rangle | \le t {\delta}/2.\] Successive
approximation of any point on the sphere by points in the net and
homogeneity of the inequality (\ref{fl2}) completes the proof, where
$C_3 = \sqrt{2\ln 5 /c}$. $\hfill \square$

We will also need a covering result of Sch\"{u}tt \cite{Schutt},
about the covering number of the unit ball of $\ell_1^m$ by
euclidean balls: There exists a universal constant $C_5$ such that
for every $k<m$
\begin{equation}\label{Schutt}
N\left(\sqrt{m} B(\ell_1^m), C_4\sqrt{{m\over k} \log {m \over
k}}B(\ell_2^m)\right) \le e^k
\end{equation}
where for two convex bodies $K$ and $T$ the
number $N(K, T)$ denotes the minimal number of translates of $T$
needed to cover $K$.

\vskip 6pt

\noindent{\em Proof of Theorem \ref{nice}} We define
\[\Sigma_{\gamma} = \{ x \in  S^{n-1}:
{1\over \sqrt{n}}\sum_{j=1}^{n+\delta n/2} |\langle x, v_j \rangle |
\le \gamma \}\] (notice that we only use $v_j$ and not $w_j$). If a
point on the sphere is not in $\Sigma_{\gamma}$ then a lower bound
$\gamma$ holds for this point for the left hand side of
(\ref{good}). We denote by $A$ the $(n+\delta n/2) \times n$ matrix
with rows $v_j$, and for convenience denote $m = n+\delta n/2$.

We now use (\ref{Schutt}) to cover ${\rm Im} A \cap
\sqrt{n}B(\ell_1^m)$ by $e^{c_5'\delta n}$ balls of radius $r =
C_4\sqrt{({1+\delta \over c_5'\delta}) \log ({1 +\delta \over
c_5'\delta})}$ (where $c_5'$ is a universal constant to be
determined later). We have used the fact that taking a section only
reduces the covering number by euclidean balls. Denote by $y_j\in
\R^m\cap {\rm Im}A$ the centers of this covering, and let $x_j\in
\R^n$ be their pre-images, so that $Ax_j = y_j$. Since there are
only $e^{c_5'\delta n}$ of them, we can use Chernoff inequality in
the following way: For a suitably chosen universal $c_5$ the
probability that for a single $i$ we have $| \langle x_j, w_i
\rangle |\ge 3c_{5}|x_j|/\sqrt{n}$ is greater than $1/2$ (this is
not difficult to prove, see for example \cite{AFM}). Therefore, by
Chernoff's theorem, the probability that for at least $1/3$ of the
indices $i= 1, \ldots, \delta n/2$ we have that $| \langle x_j, w_i
\rangle |\ge 6c_{5}|x_j|/\sqrt{n}$ is greater than $1 -
e^{-2c_5'\delta n}$ (this is our definition of $c_5'$, which is
universal). We get that with probability $1 - e^{-c_5'\delta n}$ for
every $j$ we have
\[ c_5 \delta |x_j| \le {1\over \sqrt{n}}\sum_{i=1}^{\delta n/2}|\langle x_j, w_i\rangle|.
\]

Let $x\in \sph$, and consider \begin{equation}\label{nicer} {1\over
\sqrt{n}}\sum_{i=1}^{n+\delta n/2} | \langle x, v_i \rangle | +
{1\over \sqrt{n}}\sum_{i=1}^{\delta n/2} | \langle x, w_i \rangle |
\end{equation}
(which is the same as the left hand side of (\ref{good}) up to a
factor $(1+\delta)$). Recall that if $x\in \sph$ and $x\not \in
\Sigma_{\gamma}$, we have a lower bound at least $\gamma$ for
(\ref{nicer}). Otherwise, we have $Ax \in \gamma
\sqrt{n}B(\ell_1^m)$ (and of course also $Ax \in {\rm Im}A$).
Therefore, there is some index $j$ with $|Ax - \gamma Ax_j|<\gamma
r$, where we use absolute value $|\cdot |$ to denote the euclidean
norm. This implies, using Proposition \ref{sasha} (which holds with
%here keep an eye for changes in sasha's est.
probability at least $1- e^{-c_1'\delta n^{1/6}}$), that $|x -
\gamma x_j|<\gamma r/(c_1\delta)$. In particular, $|x_j|> {1\over
\gamma} - r/(c_1\delta)$. By (\ref{fl2}) we know  that this implies
that

\begin{eqnarray*}
{1\over \sqrt{n}} \sum_{i=1}^{\delta n/2} | \langle x, w_i \rangle |
&\ge& {1\over \sqrt{n}} \sum_{i=1}^{\delta n/2} | \langle \gamma
x_j, w_i \rangle |  - {1\over \sqrt{n}} \sum_{i=1}^{\delta n/2} |
\langle
x-\gamma x_j, w_i \rangle | \\
&\ge& c_5 \delta \gamma |x_j| - C_5 \sqrt{\delta}\gamma r/(c_1\delta) \\
&\ge & c_5\delta - r\gamma (1+ C_5\sqrt{\delta})/(c_1\delta).
\end{eqnarray*}
This tells us we may choose $\gamma = \delta^2 c_5c_1/ ( 2r
(1+C_3\sqrt{\delta}))$, and have a lower bound $c_5 \delta/2$ for
this set. For the other set we have lower bound $\gamma$, that is
(remembering what was $r$), $c_0\delta^{5/2}/\log(1/\delta)^{1/2}$
(for $c_0$ a universal constant suitably chosen). The proof is
complete. $\hfill \square$

{\noindent Shiri Artstein-Avidan, Department of Mathematics,
Princeton University, Fine Hall, Washington Road, Princeton NJ
08544-1000 USA and School of Mathematics, Institute for Advanced
Study, 1 Einstein Drive, Princeton NJ 08540 USA.

\vskip 2pt

\noindent Email address: artstein@princeton.edu

\vskip 12pt

\noindent Omer Friedland, School of Mathematical Science, Tel Aviv
University, Ramat Aviv, Tel Aviv, 69978, Israel.

\vskip 2pt

\noindent Email address: omerfrie@post.tau.ac.il

\vskip 12pt

\noindent Vitali Milman, School of Mathematical Science, Tel Aviv
University, Ramat Aviv, Tel Aviv, 69978, Israel.

\vskip 2pt

\noindent Email address: milman@post.tau.ac.il

\vskip 12pt

\noindent Sasha Sodin, School of Mathematical Science, Tel Aviv
University, Ramat Aviv, Tel Aviv, 69978, Israel.

\vskip 2pt

\noindent Email address: sodinale@post.tau.ac.il}

\end{document}